\documentclass{llncs}
\newcommand{\ignore}[1]{}
\newcommand{\mod}{\ensuremath{\mathrm{mod}}}
\newcommand{\boxi}{\ensuremath{\mathrm{box}}}
\newcommand{\bbox}{\rule{0.6em}{0.6em}}
\newcounter{obs}
\newenvironment{obs}{\refstepcounter{obs}\vspace{0.05in}\par\noindent
\textit{Observation $\theobs$. }}{\par}
\newenvironment{keywords}{\par\textit{Keywords: }}{}
\newcounter{claim}
\renewenvironment{claim}{\refstepcounter{claim}\vspace{0.1in}\par\noindent
\textbf{Claim \theclaim.}}{}

\begin{document}
\title{Boxicity of Halin Graphs\thanks{The work done by the first and second
authors was partially supported by a DST grant SR/S3/EECE/62/2006.}}
\author{L. Sunil Chandran\inst{1} \and
Mathew C. Francis\inst{1} \and
Santhosh Suresh\inst{2}}
\institute{Dept. of Computer Science and Automation, Indian Institute of
Science, Bangalore--560012, email: \{sunil,mathew\}@csa.iisc.ernet.in. 
\and
Dept. of Mechanical Engineering, Indian Institute of Technology Madras,
Chennai--600036, email: santhosh.suresh@gmail.com}
\maketitle
\bibliographystyle{plain}
\begin{abstract}
A $k$-dimensional box is the Cartesian product $R_1\times R_2\times\cdots
\times R_k$ where each $R_i$ is a closed interval on the real line.
The boxicity of a graph $G$, denoted as $\boxi(G)$ is the minimum integer
$k$ such that $G$ is the intersection graph of a collection of
$k$-dimensional boxes.
Halin graphs are the graphs formed by taking a tree with no degree 2 vertex
and then connecting its leaves to form a cycle in such a way that the graph
has a planar embedding. We prove that if $G$ is a Halin graph that is not
isomorphic to $K_4$, then $\boxi(G)=2$. In fact, we prove the stronger result
that if $G$ is a planar graph formed by connecting the leaves of any tree
in a simple cycle, then $\boxi(G)=2$ unless $G$ is isomorphic to $K_4$ (in
which case its boxicity is 1).
\end{abstract}
\section{Introduction}
Let $\mathcal{F}=\{S_i~|~i\in V\}$ be a collection of subsets of a universe
$U$ where $V$ is an index set. The graph $\Lambda(\mathcal{F})=(V,E)$ where
$E=\{(i,j)~|~
S_i\cap S_j\not=\emptyset\}$ is called the \emph{intersection graph} of 
$\mathcal{F}$. When $\mathcal{F}$ is a collection of intervals on the real
line, $\Lambda(\mathcal{F})$ is called an \emph{interval graph}.

A $k$-dimensional box or $k$-box in short is the Cartesian product
$R_1\times R_2\times\cdots\times R_k$ where each $R_i$ is a closed interval
on the real line. Two $k$-boxes, $(P_1,\ldots,P_k)$ and $(Q_1,\ldots,Q_k)$
are said to have a non-empty intersection if $P_i\cap Q_i\not=\emptyset$,
for $1\leq i\leq k$. The boxicity of a graph $G$, denoted as $\boxi(G)$,
is defined to be the minimum integer $k$ such that $G$ is the intersection
graph of a collection of $k$-boxes. Since 1-boxes are nothing but closed
intervals on the real line, interval graphs are the graphs with boxicity at
most 1. We take the boxicity of a complete graph to be 1.

For a graph $G=(V,E)$, we write $G=T\cup C$ if
$E=E(T)\cup E(C)$ where
$T$ is a tree on
the vertex set $V$ and $C$ is a simple cycle on
the leaves of $T$.
$G$ is called a Halin graph if $G$ has a planar embedding and $T$ has no
vertices of degree 2.
The notion of Halin graphs were first used by Halin \cite{Halin} in his
study of minimally 3-connected graphs. Bondy and Lovasz proved that these
graphs are almost pancyclic -- they contain a cycle of each length
between 3 and $n$ with the possible exception of one length, which must be
even. Bondy has also shown that every Halin graph is 1-hamiltonian. Lovasz
and Plummer \cite{LovPlum} show that every Halin graph with an even number
of vertices is minimal bicritical (a graph is bicritical if the removal of any
two vertices from the graph will result in a graph with a perfect matching).
Halin graphs are also interesting because some problems that are NP-complete
for general graphs have been shown to be polynomial-time solvable for Halin
graphs. Examples are the travelling salesman problem \cite{Cornuejols} and the
problem of finding a dominating cycle with at most $l$ vertices
\cite{Skowronska}. 

It has been shown in \cite{TrotHar} that every
Halin graph is a 2-interval graph -- i.e., the intersection graph of sets,
each of which is the union of at most 2 intervals.
We show in this paper that the boxicity
of a Halin graph (not isomorphic to $K_4$) is equal to 2 which means that
every Halin graph is the
intersection graph of axis-parallel rectangles on the plane as well. In fact,
we show a stronger result
-- we show that our result holds for any graph $G=T\cup C$ that has a planar
embedding, even if there are vertices of degree 2 in $T$. Since $\boxi(G)=1$
when $G$ is isomorphic to $K_4$, we show our result for graphs not
isomorphic to $K_4$.

The concept of boxicity, introduced by F. S. Roberts \cite{Roberts},
finds applications in fields such as ecology and operations research.
Computing the boxicity of a graph was shown to be NP-hard by Cozzens \cite{Coz}.
This was improved by Yannakakis \cite{Yan1} and later by Kratochvil
\cite{Kratochvil} who showed that deciding whether the boxicity of a graph is
at most 2 itself is NP-complete.
An upper bound on the boxicity of general graphs is given in
\cite{CFNMaxdeg} where it is shown that $\boxi(G)\leq 2\Delta^2$ where $G$
is any graph and $\Delta$ is the maximum degree of a vertex in $G$.
Also, for any graph $G$ on $n$ vertices and maximum degree $\Delta$, 
$\boxi(G)\leq \lceil(\Delta+2)\ln n\rceil$ \cite{tech-rep}.
It was shown in \cite{CN05} that
$\boxi(G)\leq \mbox{tw}(G)+2$ where $\mbox{tw}(G)$ is the treewidth of $G$.
Upper bounds on the boxicity of some special classes of
graphs such as chordal graphs, circular-arc graphs, AT-free graphs, permutation
graphs and co-comparability graphs are also given in \cite{CN05}.
The boxicity of planar graphs was shown to be at most 3 by Thomassen
\cite{Thoma1}. A better bound holds for outerplanar graphs, a subclass of
planar graphs. Scheinerman \cite{Scheiner} showed that the boxicity of
outerplanar graphs is at most 2. But this bound does not hold for the
class of series-parallel graphs, a slightly bigger subclass of planar graphs
than the outerplanar
graphs. Bohra et al. \cite{CRB1}, showed that there exists series-parallel
graphs with boxicity 3. In this paper, we consider another subclass of planar
graphs, namely the class of Halin graphs. We show that the boxicity of Halin
graphs is at most 2.

\section{Definitions and Notations}
As mentioned before, the notation $G=T\cup C$ is used to denote a graph that
is formed by connecting the leaves of a tree $T$ so that the subgraph
induced by the leaves of $T$ in $G$ is the simple cycle $C$. $G$ is
called a Halin graph if $G$ is planar and $T$ has no vertex of degree 2.
For a graph $G$, $V(G)$ and $E(G)$ denote the vertex set and edge set of $G$
respectively. For a vertex $u\in V(G)$, $N_G(u)=\{v\in V(G)~|~ (u,v)\in E(G)
\}$. This is often abbreviated to just $N(u)$ when the graph under
consideration is clear. Given $H\subseteq V(G)$, we denote by $G_H$ the
subgraph induced by the vertices of $H$ in $G$.

A graph $G_1$ is said to be the ``supergraph'' of a graph $G_2$
if $V(G_1)=V(G_2)$ and $E(G_1)\supseteq E(G_2)$.
Also, given two graphs $G_1$ and $G_2$ on the same vertex
set $V$ we define $G_1\cap G_2$ to be the graph with vertex set $V$ and edge
set $E(G_1)\cap E(G_2)$.
As shown in \cite{Roberts}, for any graph $G$,
$\boxi(G)\leq k$ if and only if
there exists $k$ interval graphs $I_0,\ldots,I_k$ such that $G=I_0\cap\cdots
\cap I_k$. Note that for this, each $I_i$ has to be a supergraph of $G$.

\section{Our result}
We have our main result as the following theorem.
\begin{theorem}
If $G=T\cup C$, where $T$ is a tree and $C$ is a simple cycle of the leaves of
$T$ such that $G$ is planar, then $\boxi(G)=2$ if $G$ is not isomorphic to
$K_4$.
\end{theorem}
\begin{corollary}
Every Halin graph has boxicity equal to $2$ unless it is isomorphic to $K_4$,
in which case it has boxicity equal to $1$.
\end{corollary}

\section{The proof}
Let $G=T\cup C$ where $C$ is a simple cycle connecting the leaves
of a tree $T$ such that $G$ is planar. Our strategy will be to construct
two interval graphs $G_1$ and $G_2$ such that $G=G_1\cap G_2$ thus proving
that boxicity of $G$ is at most 2. We will assume that $G$ is not a wheel
since it can be seen that a wheel being just a universal vertex added to a
cycle, has boxicity 2 unless it is a $K_4$ (in which case the boxicity is
just 1).

\subsection{Finding $u'$}
Let $S=V(G)-V(C)$ denote the set of internal vertices of the tree $T$.
We claim that there is a vertex $u'\in S$ such that $|N(u')\cap S|=1$ and
$|N(u')\cap V(C)|\geq 1$.
If there is no such vertex, then $G_S$, the induced subgraph of $G$ on $S$, has
no vertices of degree 1 which is not possible since $G_S$ is a tree ($G_S$
has more than one vertex since $G$ is not a wheel).
Now, $u'$ has at least one leaf as its neighbour since if it did not, then
its degree in $T$ is 1 implying that $u'$ is a leaf of $T$ -- a contradiction.

\subsection{Fixing the root of $T$}
Designate the internal vertex of $T$ adjacent to $u'$, say $r$,
to be the root of $T$.
Given two vertices $u$ and $v$, $u$ is said to be
an \emph{ancestor} of $v$
if $u$ lies in the path $rTv$ and $u$ is said to be a \emph{descendant} of
$v$ if $v$ is an ancestor of $u$. Note that every vertex is an ancestor
and a descendant of itself.
Let $D(u)$ for any vertex $u\in V(G)$ be defined as the set of all leaves
of $T$ that are descendants of $u$. It can be easily seen that if $u$ is a
descendant of $v$, then $D(u)\subseteq D(v)$.

\subsection{Ordering the vertices of $C$}
Let $|V(C)|=k$ and
let $C$ be $p_0 p_1\ldots p_{k-1} p_0$. Note that $D(u')$ cannot contain all
the leaves since that would mean that $D(u')=D(r)$, implying that $u'$ is the
only neighbour of $r$ in $T$. Then the degree of $r$ in $T$ would be 1, a
contradiction since $r$ is an internal vertex in $T$ and not a leaf.
Therefore, we can always find a leaf $p_i\in D(u')$ such that
$p_{(i-1)\mod k}\not\in D(u')$ (recall that $u'$ has at least one leaf as its
neighbour and therefore, $D(u')$ is not empty).
We define $l_j=p_{(i+j)\mod k}$, for $0\leq j\leq k-1$.
This implies that $l_{k-1}\not\in D(u')$ since $l_{k-1}=p_{(i-1)\mod k}$.
For $u\in V(C)$, we define $c(u)=i$ when $u=l_i$.\\

\noindent For the convenience of the reader, we summarize the construction
as of now:
\begin{itemize}
\item We chose a vertex $u'$ such that its neighbourhood contains
exactly one internal vertex and at least one leaf of $T$.
\item We chose the only internal vertex in the neighbourhood of $u'$ to be the
root $r$ of $T$ and defined the natural tree-order on $T$ with $r$ as the root.
We also defined $D(u)$ to be the set of all leaves that are descendants (in
our tree-order) of the vertex $u$.
\item We defined a linear ordering $l_0,\ldots,l_{k-1}$ of the vertices in 
$V(C)$ (the leaves of $T$) where $l_0\in D(u')$ and $l_{k-1}\not\in D(u')$.
\end{itemize}
\begin{claim}\label{ourclaim}
For any vertex $u\in V(G)$, the vertices in $D(u)$ will
occur in consecutive places in the ordering $l_0,\ldots,l_{k-1}$ of the vertices
in $C$. In other words, if $u\in V(G)$ and $x,y,z\in V(C)$ such that
$c(x)<c(z)<c(y)$ then it is not possible that $x,y\in D(u)$ and $z\not\in
D(u)$.
\end{claim}

Though the statement of the claim looks intuitive, its proof involves some
technical details. Therefore, the reader might choose to skip the proof and
continue with the rest of the construction so as not to get distracted from
the main theme.

\begin{proof}[Claim \ref{ourclaim}]
If $u$ is a leaf of $T$, then the claim is trivially true.
Let us assume that this is not the case.

Consider any planar embedding of $G$. The cycle $C$ divides the plane into
a bounded region and an unbounded region. We claim that all the internal
vertices of $T$ will lie in one of these regions. Suppose there are two
internal vertices of $T$ such that they lie in different regions of $C$.
Then, the path between them in $T$ will
have to pass the boundary of $C$. But the path cannot pass through a leaf of
$T$ and because the drawing is planar,
no edge of the path can cross the boundary of $C$. We thus have a
contradiction.
Therefore, $C$ forms the boundary of a face in any planar drawing of $G$.

Now, consider a planar embedding of $G$ such that $C$ forms the boundary of
the unbounded face (i.e., all the internal vertices of $T$ lie in the
bounded region of $C$).
Suppose $x,y\in D(u)$ and $z\not\in D(u)$ such that
$c(x)<c(z)<c(y)$ (recall that $c(l_i)=i$). Let $B=xCyTuTx$. It can be easily verified that $B$ has
exactly two regions -- one bounded and the other unbounded.
We say that a vertex is
``inside'' $B$ if it lies in the bounded region bounded by $B$ and say that
it is ``outside'' $B$ if it lies in the unbounded region whose boundary is $B$.
We say that a vertex ``lies on'' $B$ if it is in $B$.

\begin{obs}\label{obsleafpos}
Because of the planar
embedding of $G$ that we have chosen, it can be seen that any leaf vertex will
have to either lie on $xCy$ or outside $B$. 
\end{obs}
\begin{obs}\label{obsrnotonB}
$r$ does not lie on $B$.
\end{obs}
We can assume that $r\not=u$ since that would contradict our assumption
that $z\not\in D(u)$. Also, $r$ cannot lie on $yTu$ or $uTx$
since it contradicts our assumption that $x$ and $y$ are descendants of $u$
and it cannot lie on
$xCy$ since it is not a leaf. Therefore, $r$ does not lie on $B$.

\begin{obs}\label{obsu'pos}
$u'$ is not inside $B$.
\end{obs}
If $u'$ is
inside $B$, then $l_0$ cannot be outside $B$ since $u'$ is adjacent to $l_0$.
From Observation \ref{obsleafpos}, $l_0$ is in $xCy$ which implies that $x=l_0$ and $u'$, being the
only internal vertex in $N(l_0)$, should lie on $uTx$.
This contradicts our assumption that $u'$ is inside $B$.

\begin{obs}\label{obsroutsideB}
$r$ is outside $B$.
\end{obs}
Now suppose $r$ is inside $B$. Then, $u'$ cannot be outside $B$ since
since $r$ is adjacent to $u'$ and it cannot be inside $B$ due to Observation
\ref{obsu'pos}.
Therefore, $u'$ lies on $B$.
If $u\not=u'$, then the fact $r$ is the only internal vertex adjacent to $u'$
implies that $r$ will have to lie on $B$, which contradicts Observation
\ref{obsrnotonB}.
Therefore, $u=u'$. Now, it can be seen that because of our choice of $u'$ and
$r$, $D(u')=N(u')-\{r\}$. This means that $uTx$ and $uTy$ are the
edges $u'x$ and $u'y$ respectively and therefore, any path from $r$ (inside
$B$) to a vertex outside $B$ will have to go through $u'$. Now, consider the
leaf $l_{k-1}$. By our construction, $l_{k-1}\not\in D(u')$. Therefore,
$y\not=l_{k-1}$ and $l_{k-1}$ does not lie on $xCy$ and hence lies outside $B$
(from Observation \ref{obsleafpos}).
The path from $r$ to $l_{k-1}$ will have to go through $u'$ as we have noted
before -- but this implies that $l_{k-1}\in D(u')$ which is
a contradiction. Therefore, $r$ is outside $B$ since we know from Observation
\ref{obsrnotonB} that
$r$ does not lie on $B$.

Because of Observation \ref{obsroutsideB}, the path $zTr$ must contain a vertex $v$ in $B$ because of 
our planarity assumption. But if $v\not=u$, then $x$ and $y$ cannot both be
descendants of $u$ since either $rTx$ or $rTy$ will not contain $u$.
If $v=u$, then $rTz$
contains $u$ and therefore, $z\in D(u)$, again a contradiction.

This proves our claim that for any vertex $u\in V(G)$, the vertices in $D(u)$
have to occur consecutively in the ordering $l_0,l_1,\ldots,l_{k-1}$.\hfill\qed
\end{proof}

\subsection{Construction of the interval graphs $G_1$ and $G_2$}
We define $f_1$ and $f_2$ to be mappings of the vertex set $V(G)$ to closed
intervals on the real line. Let $G_1$ and $G_2$ denote the interval graphs
defined by $f_1$ and $f_2$ respectively.

For a vertex $u\in V(G)$, let $d(u)$ denote the number of ancestors of $u$
other than itself (or ``depth'' of $u$ in $T$).
Let $h$ denote the maximum depth of a vertex in $T$. Recall that $k=|V(C)|$ and
$S$ denotes the set of internal vertices of $T$.\\

\noindent Definition of $f_1$:\\
For $u\in V(G)$,\\
\indent $f_1(l_0)=[0,k]$.\\
\indent $f_1(u)=[c(u)-1/2,c(u)+1/2]$, if $u\in V(C)$ and $u\not=l_0$.\\
\indent $f_1(u)=[\min_{v\in D(u)}\{c(v)\},\max_{v\in D(u)}\{c(v)\}]$, if $u
\in S$.\\

\noindent Definition of $f_2$:\\
For $u\in V(G)$,\\
\indent $f_2(u')=[d(u'),h+2]$.\\
\indent $f_2(u)=[d(u),d(u)+1]$, if $u\in S$ and $u\not=u'$.\\
\indent $f_2(l_0)=[h+2,h+2]$.\\
\indent $f_2(l_1)=[d(l_1),h+2]$.\\
\indent $f_2(l_{k-1})=[d(l_{k-1}),h+2]$.\\
\indent $f_2(u)=[d(u),h+1]$, if $u\in V(C)$ and $u$ is not $l_0$,
$l_1$ or $l_{k-1}$.\\

\begin{claim}\label{G_1supgraph}
$G_1$ is a super graph of $G$.
\end{claim}
\begin{proof}
Consider an edge $(u,v)\in E(G)$. Clearly, $(u,v)\in E(T)$ or $(u,v)\in E(C)$.
\begin{enumerate}
\item $(u,v)\in E(T)$.\\
\hspace*{0.2in}In this case, either $u$ is an ancestor of $v$ or vice versa
as $T$ is a tree. Let us assume without loss of generality that $u$ is the
ancestor of $v$. Therefore, $D(v)\subseteq D(u)$. There are two possibilities
now:
\begin{enumerate}
\item $u$ and $v$ are both internal vertices of $T$.\\
\hspace*{0.2in} Since $D(v)\subseteq D(u)$, we have $\min_{x\in D(u)}\{c(x)\}\leq
\min_{x\in D(v)}\{c(x)\}\leq$\linebreak$\max_{x\in D(v)}\{c(x)\}\leq
\max_{x\in D(u)}
\{c(x)\}$. Therefore, $f_1(u)\cap f_1(v)\not=\emptyset$, which implies that
$(u,v)\in E(G_1)$.
\item $u$ is an internal vertex of $T$ and $v$ is a leaf vertex of $T$.\\
\hspace*{0.2in}Since $v\in D(u)$, $\min_{x\in D(u)}\{c(x)\}\leq c(v)\leq\max_{x\in
D(u)}\{c(x)\}$. Thus, both $f_1(u)$ and $f_1(v)$ contain the point $c(v)$ and
therefore, $(u,v)\in E(G_1)$ (Note that $c(l_0)=0$ and thus $c(l_0)\in
f_1(l_0)$).
\end{enumerate}
\item $(u,v)\in E(C)$.\\
\hspace*{0.2in}Without loss of generality, we can assume that $u=l_i$, for some $i$,
and $v=l_{(i+1)\mod k}$. For $1\leq i\leq k-2$, $f_1(u)$ and $f_1(v)$ contain
the point $i+1/2$. If $u=l_0$ or $v=l_0$, then it is clear that $(u,v)\in
E(G_1)$, since $f_1(l_0)$ contains $f_1(u)$,$\forall u\in V(G)$.
\end{enumerate}

Therefore, $G_1$ is a supergraph of $G$.\hfill\qed
\end{proof}

\begin{claim}\label{G_2supgraph}
$G_2$ is a supergraph of $G$.
\end{claim}
\begin{proof}
Consider an edge $(u,v)\in E(G)$. We have the following three cases now.
\begin{enumerate}
\item $u$ or $v$ is $l_0$.\\
\hspace*{0.2in} By our choice of $l_0$, it is adjacent only to $l_1$, $l_{k-1}$
and $u'$. Since $f_2(l_0)$, $f_2(l_1)$, $f_2(l_{k-1})$ and $f_2(u')$ contain the
point $h+2$, all the edges incident on $l_0$ in $G$ are also present in $G_2$.
\item $(u,v)\in E(T)$, $u\not=l_0$ and $v\not=l_0$.\\
\hspace*{0.2in} Let us assume without loss of generality that $u$ is the parent
of $v$. It is easily seen that $d(v)=d(u)+1$. Since $u\not=l_0$ and $v\not=l_0$,
the point $d(u)+1$ is contained in both $f_2(u)$ and $f_2(v)$ (Recall that
$d(u)\leq h$, $\forall u\in V(G)$).
\item $(u,v)\in E(C)$, $u\not=l_0$ and $v\not=l_0$.\\
\hspace*{0.2in} Since $u$ and $v$ are leaf vertices, $f_2(u)$ and $f_2(v)$ both
contain the point $h+1$ and therefore $(u,v)\in E(G_2)$.
\end{enumerate}

This shows that $G_2$ is a supergraph of $G$.\hfill\qed
\end{proof}

\begin{claim}\label{missedges}
$G=G_1\cap G_2$.
\end{claim}
\begin{proof}
Since Claims \ref{G_1supgraph} and \ref{G_2supgraph} have established that
$G_1$ and $G_2$ are supergraphs of $G$, it is sufficient to show that for
any pair of vertices $(u,v)\not\in E(G)$, $(u,v)\not\in E(G_1)$ or
$(u,v)\not\in E(G_2)$. Consider such a pair of vertices. There are three cases
to be considered.
\begin{enumerate}
\item One of $u$ or $v$ is $l_0$.\\
\hspace*{0.2in} Let us assume without loss of generality that $u=l_0$.
$(u,v)\not\in E(G)$ now implies that $v\in V(G)-\{l_1,l_{k-1},u'\}$
since $l_0$ is only adjacent to $l_1$, $l_{k-1}$ and $u'$ in $G$.
It can be easily verified that only $f_2(u')$, $f_2(l_1)$ and
$f_2(l_{k-1})$ have a non-empty intersection with $f_2(l_0)$. Therefore,
$(u,v)\not\in E(G_2)$.
\item $u\not=l_0$, $v\not=l_0$ and one of $u$ and $v$ is the ancestor of the
other.\\
\hspace*{0.2in} Let us assume without loss of generality that $u$ is the
ancestor of
$v$. This implies that $d(v)\geq d(u)+2$ since $(u,v)\not\in E(G)$. We know that
$u\not=u'$ since all the descendants of $u'$ are its neighbours by our
choice of $u'$ and the root $r$. Now, since
$u\not=u'$, the right end-point of $f_2(u)$ is $d(u)+1$ and for all possible
choices of $v$ (excluding $l_0$), the left end-point of $f_2(v)$ is
$d(v)\geq d(u)+2$. Therefore, $f_2(u)\cap f_2(v)=\emptyset$ by the definition
of $f_2$. Thus, in this case, $(u,v)\not\in E(G_2)$.
\item $u\not=l_0$, $v\not=l_0$ and neither one of $u$ and $v$ is an ancestor
of the other.\\
\hspace*{0.2in} One of the following three subcases hold.
\begin{enumerate}
\item $u$ and $v$ are both leaves of $T$.\\
\hspace*{0.2in} Let $u=l_i$ and $v=l_j$. Assume without loss of generality that
$i<j$. Since neither of $u$ or $v$ is $l_0$, we have $1\leq i<j\leq k-1$. 
Also, $j>i+1$ as $(u,v)\not\in E(G)$. Therefore, $f_1(l_i)\cap f_1(l_j)=
\emptyset$, from the definition of $f_1$. Thus, we have $(u,v)\not\in E(G_1)$.
\item $u$ and $v$ are both internal vertices of $T$.\\
\hspace*{0.2in} Since $u\not\in rTv$ and $v\not\in rTu$, we have $D(u)\cap D(v)
=\emptyset$ (To see this, suppose there is a vertex $z\in D(u)\cap D(v)$. Then
both $u$ and $v$ would lie on $rTz$, implying that either $u\in rTv$ or $v\in
rTu$). Now, from Claim \ref{ourclaim}, we have $\max_{x\in D(u)}\{c(x)\}<
\min_{x\in D(v)}\{c(x)\}$ or $\max_{x\in D(v)}\{c(x)\}<\min_{x\in D(u)}
\{c(x)\}$. By the definition of $f_1$,
it can be seen that $f_1(u)\cap f_1(v)=\emptyset$, implying that $(u,v)
\not\in E(G_1)$.
\item One of $u$ and $v$ is a leaf of $T$ and the other is an internal vertex
of $T$.\\
\hspace*{0.2in} Let us assume that $u$ is an internal vertex and $v$ is a leaf
of $T$. Since we are considering the case when neither of $u$ and $v$
is an ancestor of the other and neither is $l_0$, we have $v\not\in D(u)$ and
$v\not=l_0$. From Claim \ref{ourclaim}, we know that
either $c(v)<\min_{x\in D(u)}\{c(x)\}$ or $c(v)>\max_{x\in D(u)}\{c(x)\}$.
Therefore, by definition of $f_1$ and because $v\not=l_0$,
$f_1(u)\cap f_1(v)=\emptyset$ and thus we have $(u,v)\not\in E(G_1)$.
\end{enumerate}
\end{enumerate}

Since we have considered all possible cases when $(u,v)\not\in E(G)$ and have
shown that in each case, $(u,v)$ is not present either in $E(G_1)$ or in
$E(G_2)$, it follows that $G=G_1\cap G_2$.\hfill\qed

Now, to complete the proof, we show that if $G$ is not
isomorphic to $K_4$, then $\boxi(G)\geq 2$.
Suppose $G$ is not isomorphic to $K_4$.
We will show that $G$ is not an interval graph.
By definition of $G$, $|V(C)|\geq 3$. If $|V(C)|>3$, then $C$ is an
induced cycle with more than 3 vertices which means that $G$ cannot be an
interval graph and therefore $\boxi(G)\geq 2$.
If $|V(C)|=3$, then $C$ is a triangle.
Now, all the leaves in $V(C)$ cannot be adjacent to the same internal vertex of
$T$. To see this, look at $G_S$, the subgraph induced by $S$ in $G$ (recall
that $S=V(G)-V(C)$, or the set of internal vertices of $T$).
Since $G$ is not isomorphic to $K_4$, $G_S$ is a tree with more than one vertex.
Therefore, there are at least two vertices of degree 1 if $G_S$.
But since all the vertices in $V(C)$ are adjacent
only to one vertex of $S$ in $G$, there should be at least one vertex in $S$
with degree 1 in $G$ -- which is a contradiction since all vertices
of $S$, being internal vertices of $T$, have degree more than 1 in $G$.
Therefore, we can find two leaves, say $x$ and $y$, of $T$ such that
they are adjacent to different internal vertices in $T$. Let $u$ and $v$ denote
the internal vertices of $T$ adjacent to $x$ and $y$ respectively.
Now, $xuTvyx$ forms an induced cycle of length greater than or equal to 4.
Therefore, $G$ cannot be an interval graph. Thus, we have $\boxi(G)\geq 2$.
\hfill\bbox
\end{proof}

\end{document}